\documentclass[12pt]{amsart}
\usepackage[utf8]{inputenc}
\usepackage{amssymb, amsmath}
\usepackage{fourier}
\usepackage{enumerate}

\theoremstyle{plain}
\newtheorem*{Theorem}{Theorem}

\newtheorem*{Proposition}{Proposition}

\theoremstyle{remark}
\newtheorem*{Remark}{Remark}

\newtheorem*{Example}{Example}

\newcommand{\bit}{\begin{itemize}}
\newcommand{\eit}{\end{itemize}}
\newcommand{\ben}{\begin{enumerate}}
\newcommand{\een}{\end{enumerate}}
\newcommand{\be}{\begin{equation}}
\newcommand{\ee}{\end{equation}}
\newcommand{\ba}{\begin{array}}
\newcommand{\ea}{\end{array}}

\newcommand{\abs}[1]{\left|#1\right|}
\newcommand{\norm}[1]{\left|\left|#1\right|\right|}

\newcommand{\dt}{\mathrm{d}t}

\newcommand{\eps}{\varepsilon}

\newcommand{\inner}[2]{\left\langle #1,#2 \right\rangle}

\newcommand\cA{\mathcal A}
\newcommand\cB{\mathcal B}

\newcommand\cH{\mathcal H}
\newcommand\cK{\mathcal K}

\newcommand\C{\mathbb C}
\newcommand\R{\mathbb R}
\newcommand\N{\mathbb N}

\newcommand\Det{\operatorname{Det}}
\newcommand{\ler}[1]{\left( #1 \right)}

\title[Centrality and local monotonicity]{Connections between centrality and local monotonicity of certain functions on $C^*$-algebras}

\author{D\'aniel Virosztek}
\address{Department of Analysis, Institute of Mathematics\\
Budapest University of Technology and Economics\\
H-1521 Budapest, Hungary
and
MTA-DE ``Lend\" ulet'' Functional Analysis Research Group, Institute of Mathematics\\
         University of Debrecen\\
         H-4002 Debrecen, P.O. Box 400, Hungary}
\email{virosz@math.bme.hu}
\urladdr{http://www.math.bme.hu/\~{}virosz}

\thanks{
The author was supported by the ``Lend\" ulet'' Program (LP2012-46/2012) of the Hungarian Academy of Sciences and by the National Research, Development and Innovation Office -- NKFIH, Grant No. K104206.}


\begin{document}

\begin{abstract}
We introduce a quite large class of functions (including the exponential function and the power functions with exponent greater than one), and show that for any element $f$ of this function class, a self-adjoint element $a$ of a $C^*$-algebra is central if and only if $a \leq b$ implies $f(a) \leq f(b).$ That is, we characterize centrality by local monotonicity of certain functions on $C^*$-algebras. Numerous former results (including works of \emph{Ogasawara, Pedersen, Wu,} and \emph{Moln\'ar}) are apparent consequences of our result.
\end{abstract}
\maketitle

\section{Introduction}
Connections between the commutativity of a $C^*$-algebra $\cA$ and the monotonicity of some functions defined on some subsets of $\cA$ have been investigated widely. The first result related to this topic is due to \emph{Ogasawara} who showed in 1955 that a $C^*$-algebra $\cA$ is commutative if and only if the square function is monotone on the positive cone of $\cA$ \cite{oga}. It was observed later by \emph{Pedersen} that the above statement remains true for any power function with exponent greater than one \cite{pedersen}. \emph{Wu} proved a similar result for the exponential function in 2001 \cite{wu}. Ji and Tomiyama showed in 2003 that for any function $f$ which is monotone but not matrix monotone of order $2,$ a $C^*$-algebra $\cA$ is commutative if and only if $f$ is monotone on the positive cone of $\cA$ \cite{ji-tom}. The reader is advised to consult the papers \cite{proc-eston} and \cite{nov-tik} for other closely related results.
\par
Very recently, Moln\'ar proved a local theorem, namely, that a self-adjoint element $a$ of a $C^*$-algebra $\cA$ is central if and only if $a  \leq b$ implies $\exp{a} \leq \exp{b}$ \cite{mol-char}.
\par
Motivated by the work of Moln\'ar, we show the following.
If $I=(\gamma, \infty)$ is a real interval and $f$ is a continuously differentiable function on $I$ such that the derivative of $f$ is positive, strictly monotone increasing and logarithmically concave, then a self-adjoint element $a$ of a $C^*$-algebra $\cA$ with spectrum in $I$ is central if and only if $ a \leq b$ implies $f(a) \leq f(b),$ that is, $f$ is locally monotone at the point $a.$ This result easily implies the results of Ogasawara, Pedersen, Wu, and Moln\'ar.

\section{The main theorem}
The precise formulation of our main result reads as follows (here and throughout, the symbol $\cA_s$ stands for the set of the self-adjoint elements of a $C^*$-algebra $\cA$).
\begin{Theorem} \label{fo}
Let $I=(\gamma,\infty)$ for some $\gamma \in \R\cup\{-\infty\}$ and let $f \in C^1(I)$ be such that
\ben [(i)]
\item $f'(x)>0 \qquad \ler{x \in I},$ \label{diff-poz}
\item $x<y \Rightarrow f'(x) < f'(y) \qquad \ler{x,y \in I},$ \label{diff-szig-mon}
\item $ \log \ler{f'\ler{tx+(1-t)y}} \geq t \log{f'(x)}+(1-t) \log{f'(y)} \qquad \ler{x,y \in I, \, t \in [0,1]}.$ \label{diff-log-konk}
\een
Let $\cA$ be a unital $C^*$-algebra and let $a \in \cA$ be a self-adjoint element with $\sigma(a) \subset I.$ The followings are equivalent.
\ben [(1)]
\item $a$ is central, that is, $ab=ba \qquad \ler{b \in \cA},$ \label{cent}
\item $f$ is locally monotone at the point $a,$ that is, $a \leq b \Rightarrow f(a) \leq f(b) \quad \ler{b \in \cA_s}.$ \label{lokmon}
\een
\end{Theorem}
\begin{Example}
We enumerate the most important examples of intervals and functions satisfying the conditions given in the Theorem:
\bit
\item $I=(0, \infty), \, f(x)=x^p \qquad \ler{p>1},$ \label{pl-hatv}
\item $I=(-\infty,\infty), \, f(x)=e^x.$
\eit
\end{Example}

\section{The proof of the theorem}
\paragraph*{{\bf Notation.}}
If $\varphi$ and $\psi$ are elements of some Hilbert space $\cH,$ then the symbol $\varphi \otimes \psi$ denotes the linear map $\cH \ni \xi \mapsto \inner{\xi}{\psi}\varphi \in \cH.$
\par
\medskip
The following proposition is a key step of the proof.
\begin{Proposition}
Suppose that $I=(\gamma,\infty)$ for some $\gamma \in \R\cup\{-\infty\}$ and $f \in C^1(I)$ satisfies the conditions \eqref{diff-poz}, \eqref{diff-szig-mon} and \eqref{diff-log-konk} given in the Theorem. Let $\cK$ be a two-dimensional Hilbert space, let $\{u,v\} \subset \cK$ be an orthonormal basis. Let $x,y \in I$ and set $A:= x u \otimes u + y v \otimes v.$ The followings are equivalent.
\ben[(I)]
\item $x \neq y,$ \label{nemegy}
\item there exist $\lambda, \mu \in \C$ with $\abs{\lambda}^2+\abs{\mu}^2=1$ and $t_0>0$ such that using the notation $B= \ler{u+v} \otimes \ler{u+v}$ and $w=\lambda u +\mu v$ we have
$$
\inner{f(A)w}{w}-\inner{f\ler{A+t_0 B} w}{w}>0.
$$ \label{nemmon}
\een
\end{Proposition}

\paragraph*{{\bf Notation.}}
For any fixed interval $I =(\gamma, \infty)$ and function $f \in C^1(I)$ with the properties \eqref{diff-poz}, \eqref{diff-szig-mon} and \eqref{diff-log-konk}, and different numbers $x,y \in I,$ the above Proposition provides a positive number $\inner{f(A)w}{w}-\inner{f\ler{A+t_0 B} w}{w}.$
Let us introduce
$$
\delta:=\inner{f(A)w}{w}-\inner{f\ler{A+t_0 B} w}{w}.
$$

\begin{proof}[Proof of the Proposition]
The direction \eqref{nemmon} $\Rightarrow$ \eqref{nemegy} is easy to see (by contraposition). To verify the direction \eqref{nemegy} $\Rightarrow$ \eqref{nemmon} we recall the following useful formula for the derivative of a matrix function (see \cite[Thm. 3.25]{HP_book} and also \cite[Thm. 3.33]{HP_book}).
If $A= x u \otimes u + y v \otimes v,$ then for any self-adjoint $C \in \cB(\cK)$ we have
$$
\lim_{t \to 0}\frac{1}{t} \ler{f\ler{A+tC}-f(A)}=
f'(x) \inner{Cu}{u} u \otimes u+
$$
$$
 +\frac{f(x)-f(y)}{x-y} \inner{Cv}{u} u \otimes v +\frac{f(y)-f(x)}{y-x} \inner{Cu}{v} v \otimes u+ f'(y) \inner{Cv}{v} v \otimes v.
$$
This means that for $B=\ler{u+v} \otimes \ler{u+v}$ we have

\be \label{limdef}
L:=\lim_{t \to 0}\frac{1}{t} \ler{f\ler{A+tB}-f(A)}
\ee
$$
=f'(x) u \otimes u+\frac{f(x)-f(y)}{x-y} u \otimes v +\frac{f(y)-f(x)}{y-x} v \otimes u+ f'(y) v \otimes v.
$$

The determinant of the matrix
$$
[L]=\left[\ba{cc} f'(x) & \frac{f(x)-f(y)}{x-y} \\ \frac{f(y)-f(x)}{y-x} & f'(y) \ea\right]
$$
is negative as
$
\Det [L]<0 \Leftrightarrow f'(x) f'(y)< \ler{\frac{f(x)-f(y)}{x-y}}^2 \Leftrightarrow
$
$$\Leftrightarrow \log{f'(x)}+\log{f'(y)}< 2 \log{\ler{\int_{0}^1 f'\ler{tx+(1-t)y}\dt} }.
$$
This latter inequality is true as
$$
\log{f'(x)}+\log{f'(y)}= 2 \cdot \int_{0}^1 t \log{f'(x)}+(1-t)\log{f'(y)}\dt
$$
$$
\leq 2 \int_{0}^1 \log{\ler{f'\ler{tx+(1-t)y}}}\dt < 2 \log{\ler{\int_{0}^1 f'\ler{tx+(1-t)y}\dt} }.
$$
In the above computation, the first inequality holds because of the log-concavity of $f'$ and the second (strict) inequality holds because the logarithm function is strictly concave and $f'$ is strictly monotone increasing.
\par
So, the operator $L$ (defined in eq. \eqref{limdef}) has a negative eigenvalue, that is, there exist $\lambda, \mu \in \C$ (with $\abs{\lambda}^2+\abs{\mu}^2=1$) such that with $w=\lambda u + \mu v$ we have
$$
\inner{Lw}{w}=\inner{\lim_{t \to 0}\frac{1}{t} \ler{f\ler{A+tB}-f(A)} w}{w}<0.
$$
Therefore,
$$
\lim_{t \to 0}\frac{1}{t}\ler{\inner{ f\ler{A+tB} w}{w}-\inner{ f(A) w}{w}}<0,
$$
and so there exists some $t_0>0$ such that $0<\inner{f(A)w}{w}-\inner{f\ler{A+t_0 B} w}{w}.$
\end{proof}

\begin{proof}[Proof of the Theorem]
The direction \eqref{cent} $\Rightarrow$ \eqref{lokmon} is easy to verify, because if $f$ is continuous and monotone increasing as a function of one real variable, then the map $x \mapsto f(x)$ preserves the order of commuting self-adjoint elements of a unital $C^*$-algebra.
\par
To see the contrary, assume that $a \in \cA_s, \, \sigma(a) \subset I$ and $aa'-a'a \neq 0$ for some $a' \in \cA.$ Then, by \cite[10.2.4. Corollary]{kad-ring-2}, there exists an irreducible representation $\pi: \cA \rightarrow \cB(\cH)$ such that $\pi\ler{aa'-a'a} \neq 0,$ that is, $\pi(a) \pi\ler{a'}\neq \pi\ler{a'}\pi(a).$ Let us fix this irreducible representation $\pi.$ So, $\pi(a)$ is a non-central self-adjoint (and hence normal) element of $\cB(\cH)$ with $\sigma\ler{\pi(a)}\subset I$ (as a representation do not increase the spectrum). By the non-centrality, $\sigma\ler{\pi(a)}$ has at least two elements, and by the normality, every element of $\sigma\ler{\pi(a)}$ is an approximate eigenvalue \cite[3.2.13. Lemma]{kad-ring-1}. Let $x$ and $y$ be two different elements of $\sigma\ler{\pi(a)},$ and let $\{u_n\}_{n \in \N} \subset \cH$ and $\{v_n\}_{n \in \N} \subset \cH$ satisfy
$$
\lim_{n \to \infty} \pi(a) u_n - x u_n =0, \,
\lim_{n \to \infty} \pi(a) v_n - y v_n =0, \text{ and }
\inner{u_m}{v_n}=0 \quad m,n \in \N.
$$
(As $x\neq y,$ the approximate eigenvetors can be chosen to be orthogonal.) Set $\cK_n:=\mathrm{span}\{u_n, v_n\}$ and let $E_n$ be the orthoprojection onto the closed subspace $\cK_n^{\perp} \subset \cH.$ Let
$$
\psi_n(a):=x u_n \otimes u_n + y v_n \otimes v_n + E_n \pi(a)E_n.
$$
We intend to show that
$$ 
\lim_{n \to \infty} \psi_n(a)=\pi(a)
$$
in the operator norm topology.
Let $h$ be an arbitrary non-zero element of $\cH$ and consider the orthogonal decompositions $h=h_1^{(n)}+h_2^{(n)},$ where $h_1^{(n)} \in \cK_n$ and $h_2^{(n)} \in \cK_n^{\perp}$ for any $n \in \N.$
Let us introduce the symbols $\eps_{u,n}:=\pi(a) u_n-x u_n$ and $\eps_{v,n}:=\pi(a) v_n-y v_n$ and recall that $\lim_{n \to \infty} \eps_{u,n}=0$ and $\lim_{n \to \infty} \eps_{v,n}=0$ in the standard topology of the Hilbert space $\cH.$
Now,
$$
\frac{1}{\norm{h}}
\norm{\ler{\pi(a)-\psi_n(a)}h}
$$
$$
\leq
\frac{1}{\norm{h}}
\norm{\ler{\pi(a)-\psi_n(a)}h_1^{(n)}}
+\frac{1}{\norm{h}}\norm{\ler{\pi(a)-\psi_n(a)}h_2^{(n)}}.
$$
Both the first and the second term of the right hand side of the above inequality are bounded by the term $\norm{\eps_{u,n}}+\norm{\eps_{v,n}}$ because
$$
 \frac{1}{\norm{h}}\norm{\ler{\pi(a)-\psi_n(a)}h_1^{(n)}}
 =\frac{1}{\norm{h}}\norm{\ler{\pi(a)-\psi_n(a)}\ler{\alpha_n u_n +\beta_n v_n}}
$$
$$
 =\frac{1}{\norm{h}}
\norm{\alpha_n x u_n + \alpha_n \eps_{u,n}-x \alpha_n u_n + \beta_n y v_n + \beta_n \eps_{v,n}-y \beta_n v_n}
$$
$$
\leq \frac{\abs{\alpha_n}}{\norm{h}}
\norm{\eps_{u,n}}+\frac{\abs{\beta_n}}{\norm{h}}\norm{\eps_{v,n}}
\leq \norm{\eps_{u,n}}+\norm{\eps_{v,n}}
$$
as the sequences $\left\{\abs{\alpha_n}\right\}$ and $\left\{\abs{\beta_n}\right\}$ are obviously bounded by $\norm{h},$
and
$$
\frac{1}{\norm{h}}
\norm{\ler{\pi(a)-\psi_n(a)}h_2^{(n)}}=
\frac{1}{\norm{h}}
\norm{\ler{I_\cH-E_n}\pi(a)h_2^{(n)}}
$$
$$
=\frac{1}{\norm{h}}
\norm{\ler{u_n \otimes u_n + v_n \otimes v_n}\pi(a)h_2^{(n)}}
=\frac{1}{\norm{h}}
\norm{\inner{\pi(a)h_2^{(n)}}{u_n}u_n+\inner{\pi(a)h_2^{(n)}}{v_n}v_n}
$$
$$
=
\frac{1}{\norm{h}}
\norm{\inner{h_2^{(n)}}{\pi(a)u_n}u_n+\inner{h_2^{(n)}}{\pi(a) v_n}v_n}
$$
$$
\leq
\frac{1}{\norm{h}}
\abs{\inner{h_2^{(n)}}{x u_n +\eps_{u,n}}}+ \frac{1}{\norm{h}}
\abs{\inner{h_2^{(n)}}{y v_n +\eps_{v,n}}}
$$
$$
=
\frac{1}{\norm{h}}
\abs{\inner{h_2^{(n)}}{\eps_{u,n}}}
+\frac{1}{\norm{h}}
\abs{\inner{h_2^{(n)}}{\eps_{v,n}}}
$$
$$
\leq
\frac{\norm{h_2^{(n)}}}{\norm{h}}\norm{\eps_{u,n}}
+\frac{\norm{h_2^{(n)}}}{\norm{h}}\norm{\eps_{v,n}}
\leq \norm{\eps_{u,n}}+\norm{\eps_{v,n}}.
$$
We used that $a$ is self-adjoint, hence so is $\pi(a).$
\par
So, we found that
$$
\mathrm{sup}
\left\{
\frac{1}{\norm{h}}
\norm{\ler{\pi(a)-\psi_n(a)}h}
\middle|
h \in \cH\setminus \{0\} \right\}
\leq 2 \ler{\norm{\eps_{u,n}}+\norm{\eps_{v,n}}} \to 0,
$$
which means that $\psi_n(a)$ tends to $\pi(a)$ in the operator norm topology.

\par
We have fixed $I, f, x$ and $y.$ By the Proposition, we have $\lambda, \mu \in \C$ (with $\abs{\lambda}^2+\abs{\mu}^2=1$) and $t_0>0$ such that using the notation $B_n:=(u_n+v_n) \otimes (u_n+v_n)$ and $w_n:=\lambda u_n +\mu v_n,$ we have
\be \label{deltas}
\inner{f\ler{\psi_n(a)}w_n}{w_n}-\inner{f\ler{\psi_n(a)+t_0 B_n}w_n}{w_n}=\delta>0
\ee
for any $n \in \N.$ That is, the left hand side of \eqref{deltas} is independent of $n.$
\par
The operator $B_n$ is a self-adjoint element of $\cB(\cH)$ and $\cK_n$ is a finite dimensional subspace of $\cH,$ hence by Kadison's transitivity theorem \cite[10.2.1. Theorem]{kad-ring-2}, there exists a self-adjoint $b_n \in \cA$ such that
$$
\pi\ler{b_n}_{|\cK_n}=B_{n |\cK_n}.
$$
Observe that $B_n \cK_n \subseteq \cK_n$ and so $\pi\ler{b_n} \cK_n \subseteq \cK_n.$
On the other hand, $\pi\ler{b_n}$ is self-adjoint as $b_n$ is self-adjoint, hence it follows that $\pi\ler{b_n} \cK_n^{\perp} \subseteq \cK_n^{\perp}.$ Therefore, the fact $B_n=\frac{1}{2} B_n^2$ implies that
$$
\pi \ler{\frac{1}{2}b_n^2}_{|\cK_n}
=\ler{\frac{1}{2}\pi\ler{b_n}^2}_{|\cK_n}
=\frac{1}{2} \ler{\pi\ler{b_n}_{|\cK_n}}^2
=\frac{1}{2} B_{n |\cK_n}^2=B_{n |\cK_n}.
$$
So, we can rewrite \eqref{deltas} as
\be \label{uj_delt}
\inner{f\ler{\psi_n(a)}w_n}{w_n}-\inner{f\ler{\psi_n(a)+t_0 \pi \ler{\frac{1}{2}b_n^2}}w_n}{w_n}=\delta>0
\ee
A standard continuity argument which is based on the fact that $\psi_n(a)$ tends to $\pi(a)$ in the operator norm topology shows that
\be \label{fkonv1}
\lim_{n \to \infty} \norm{f\ler{\psi_n(a)}-f\ler{\pi(a)}}=0.
\ee
Moreover, by Kadison's transitivity theorem, the sequence $\pi \ler{\frac{1}{2}b_n^2}$ is bounded (for details, the reader should consult the proof of \cite[5.4.3. Theorem]{kad-ring-1}), and hence
\be \label{fkonv2}
\lim_{n \to \infty} \norm{f\ler{\psi_n(a)+t_0\pi \ler{\frac{1}{2}b_n^2}}-f\ler{\pi(a)+t_0\pi \ler{\frac{1}{2}b_n^2}}}
=0
\ee
also holds.
By \eqref{fkonv1} and \eqref{fkonv2}, for any $\delta>0$ one can find $n_0 \in \N$ such that for $n>n_0$ we have
$$
\norm{f\ler{\psi_n(a)}-f\ler{\pi(a)}}<\frac{1}{4} \delta
$$
and
$$
\norm{f\ler{\psi_n(a)+t_0\pi \ler{\frac{1}{2}b_n^2}}-f\ler{\pi\ler{a+\frac{t_0}{2} b_n^2}}} <\frac{1}{4} \delta.
$$
Therefore, by \eqref{uj_delt}, for $n>n_0,$ the inequality
$$ 
\inner{f\ler{\pi(a)}w_n}{w_n}-\inner{f\ler{\pi\ler{a+\frac{t_0}{2} b_n^2}}w_n}{w_n}>\frac{1}{2}\delta>0
$$
holds.
In other words,
$$
f\ler{\pi(a)} \nleq f\ler{\pi\ler{a+\frac{t_0}{2} b_n^2}},
$$
or equivalently,
$$
\pi\ler{f(a)} \nleq \pi\ler{f\ler{a+\frac{t_0}{2} b_n^2}}.
$$
Any representation of a $C^*$-algebra preserves the semidefinite order, hence this means that
$$
f(a) \nleq f\ler{a+\frac{t_0}{2} b_n^2},
$$
despite the fact that $a \leq a+\frac{t_0}{2} b_n^2.$ The proof is done.
\end{proof}
\begin{Remark}
Note that our theorem generalizes Moln\'ar's result, and --- as every "local" theorem easily implies its "global" counterpart --- we recover the theorems of Ogasawara, Pedersen, and Wu, as well.
\end{Remark}
{\bf Acknowledgements.}

The author is grateful to Lajos Moln\'ar for
\bit
\item proposing the problem discussed in this paper and for sharing the manuscript version of \cite{mol-char},
\item correcting an earlier version of this paper,
\item useful discussions.
\eit
\par
The author is also grateful to J\'ozsef Pitrik for great discussions on the topic.


\begin{thebibliography}{99}

\bibitem{HP_book} F. Hiai and D. Petz, \emph{Introduction to Matrix Analysis and Applications,} Hindustan Book Agency and Springer Verlag (2014)

\bibitem{ji-tom} G. Ji and J. Tomiyama, \emph{On characterizations of commutativity of $C^*$-algebras,} Proc. Amer. Math. Soc. {\bf 131} (2003), 3845--3849.

\bibitem{kad-ring-1} R. V. Kadison and J. R. Ringrose, \emph{Fundamentals of the Theory of Operator Algebras,} Volume I, Academic Press, Orlando, 1983.

\bibitem{kad-ring-2} R. V. Kadison and J. R. Ringrose, \emph{Fundamentals of the Theory of Operator Algebras,} Volume II, Academic Press, Orlando, 1986.


\bibitem{mol-char} L. Moln\'ar, \emph{A characterization of central elements in $C^*$-algebras,} Bull. Austral. Math. Soc. {\bf 95} (2017), 138--143.

\bibitem{nov-tik} A. A. Novikov and O. E. Tikhonov, \emph{Characterization of central elements of operator algebras by inequalities,} Lobachevskii J. Math. {\bf 36} (2015), 208--210.

\bibitem{oga} T. Ogasawara, \emph{A theorem on operator algebras,} J. Sci. Hiroshima Univ. Ser. A. {\bf 18} (1955), 307--309.


\bibitem{pedersen} G.K. Pedersen, \emph{ $C^*$-Algebras and Their Automorphism Groups,} London Mathematical Society Monographs, 14, Academic Press, Inc., London-New York, 1979.

\bibitem{proc-eston} S. Silvestrov, H. Osaka and J. Tomiyama, \emph{Operator convex functions over $C^*$-algebras,} Proc. Eston. Acad. Sci. {\bf 59} (2010), 48--52.

\bibitem{wu} W. Wu, \emph{An order characterization of commutativity for $C^*$-algebras,} Proc. Amer. Math. Soc. {\bf 129} (2001), 983--987.

\end{thebibliography}
\end{document}